\documentclass{article}

\usepackage[preprint]{neurips_2025}

\usepackage[utf8]{inputenc} 
\usepackage[T1]{fontenc}    
\usepackage{hyperref}       
\usepackage{url}            
\usepackage{booktabs}       
\usepackage{amsfonts}       
\usepackage{nicefrac}       
\usepackage{microtype}      
\usepackage{xcolor}         
\usepackage{glossaries}
\usepackage{enumitem}
\usepackage{listings}
\usepackage{verbatim}
\usepackage{multirow}
\usepackage{cleveref}
\usepackage{graphicx}

\crefname{equation}{}{}
\usepackage{listings}
\newacronym{pde}{PDE}{partial differential equation}
\newacronym{ad}{AD}{automatic differentiation}
\newacronym{alcf}{ALCF}{Argonne Leadership Computing Facility}
\newacronym{blas}{BLAS}{basic linear algebra subprograms}
\newacronym{iap}{IAP}{independent activities period}
\newacronym{der}{DER}{distributed energy resource}
\newacronym{derms}{DERMS}{distributed energy resource management system}
\newacronym{minlp}{MINLP}{mixed-integer nonlinear programming}
\newacronym{nlp}{NLP}{nonlinear programming}
\newacronym{kkt}{KKT}{Karush-Kuhn-Tucker}
\newacronym{sqp}{SQP}{sequential quadratic programming}
\newacronym{ipm}{IPM}{interior-point method}
\newacronym{cpu}{CPU}{central processing units}
\newacronym{gpu}{GPU}{graphics processing units}
\newacronym{mpc}{MPC}{model predictive control}
\newacronym{ac}{AC}{alternating current}
\newacronym{dc}{DC}{direct current}
\newacronym{opf}{OPF}{optimal power flow}
\newacronym{hpc}{HPC}{high-performance computing}
\newacronym{pcg}{PCG}{projected conjugate gradient}
\newacronym{alm}{ALM}{augmented Lagrangian method}
\newacronym{dac}{DAC}{direct air capture}
\newacronym{pem}{PEM}{proton exchange membrane}
\newacronym{tea}{TEA}{technoeconomic analysis}
\newacronym{lca}{LCA}{life cycle assessment}
\newacronym{ft}{FT}{Fischer-Tropsch}
\newacronym{bess}{BESS}{battery energy storage system}
\newacronym{cqa}{CQA}{critical quality attribute}
\newacronym{lnp}{LNP}{lipid nanoparticle}
\newacronym{iso}{ISO}{independent system operator}
\newacronym{simd}{SIMD}{single instruction, multiple data}
\newacronym{mimd}{MIMD}{multiple instruction, multiple data}
\newacronym{orcd}{ORCD}{Office of Research Computing and Data}
\newacronym{mitei}{MITei}{MIT Energy Initiative}
\newacronym{hipsat}{HIP-SAT}{High School Introduction to the Physical Sciences and Advanced Technologies}
\newacronym{mpi}{MPI}{message passing interface}
\newacronym{goc}{GOC}{grid optimization competition}
\newacronym{sc}{SC}{security-constrained}
\newacronym{mp}{MP}{multi-period}
\newacronym{admm}{ADMM}{alternating direction method of multipliers}
\newacronym{gmres}{GMRES}{generalized mean residual}
\newacronym{pdhg}{PDHG}{primal-dual hybrid gradient}
\newacronym{lp}{LP}{linear programming}
\newacronym{ams}{AMS}{algebraic modeling system}
\newacronym{ncl}{NCL}{nonlinear constrained Lagrangian}
\newacronym{bcl}{BCL}{bound constrained Lagrangian}
\newacronym{sqd}{SQD}{symmetric quasi-definite}
\newacronym{spd}{SPD}{symmetric positive definite}

\title{GPU Implementation of Second-Order Linear and Nonlinear Programming Solvers}

\author{%
  Alexis Montoison\\
  Mathematics and Computer Science Division\\
  Argonne National Laboratory\\
  Lemont, IL 60439\\
  \texttt{amontoison@anl.gov}\\
  \And
  Fran\c{c}ois Pacaud\\
  Centre Automatique et Systèmes\\
  Mines Paris-PSL\\
  Paris, 75006 \\
  \texttt{francois.pacaud@minesparis.psl.eu}\\
  \And
  Sungho Shin\\
  Department of Chemical Engineering\\
  Massachusetts Institute of Technology\\
  Cambridge, MA 02139\\
  \texttt{sushin@mit.edu}\\
  \And
  Mihai Anitescu\\
  Mathematics and Computer Science Division\\
  Argonne National Laboratory\\
  Lemont, IL 60439\\
  \texttt{anitescu@mcs.anl.gov}\\
}

\begin{document}

\maketitle

\begin{abstract}
In recent years, GPU-accelerated optimization solvers based on second-order methods (e.g., interior-point methods) have gained momentum with the advent of mature and efficient GPU-accelerated direct sparse linear solvers, such as cuDSS. This paper provides an overview of the state of the art in GPU-based second-order solvers, focusing on \emph{pivoting-free interior-point methods} for large and sparse linear and nonlinear programs. We begin by highlighting the capabilities and limitations of the currently available GPU-accelerated sparse linear solvers. Next, we discuss different formulations of the Karush-Kuhn-Tucker systems for second-order methods and evaluate their suitability for pivoting-free GPU implementations. We also discuss strategies for computing sparse Jacobians and Hessians on GPUs for nonlinear programming. Finally, we present numerical experiments demonstrating the scalability of GPU-based optimization solvers. We observe speedups often exceeding 10× compared to comparable CPU implementations on large-scale instances when solved up to medium precision. Additionally, we examine the current limitations of existing approaches.

\end{abstract}

\section{Introduction}\label{eqn:intro}

This paper focuses on the implementation of solvers for problems of the following form:
\begin{align}\label{eqn:opt}
  \min_{x } \; f(x) \quad \text{s.t.} \quad g(x) \geq 0,
\end{align}
where \(x \in \mathbb{R}^n\) is the decision variable, and \(f: \mathbb{R}^n \to \mathbb{R}\) and \(g: \mathbb{R}^n \to \mathbb{R}^m\) are the smooth objective and constraint functions, respectively.
There is no loss of generality in considering only inequality constraints: equality constraints can be represented by pairs of inequalities.
We will discuss both \gls*{lp} (where \(f\) and \(g\) are affine) and \gls*{nlp} (where \(f\) and \(g\) are nonlinear), with an emphasis on algorithms designed for large, sparse instances.

Despite advances in general-purpose GPU computing, state-of-the-art mathematical programming solvers have not widely adopted these techniques. GPUs excel in repetitive computations on large data sets, such as dense matrix multiplication in AI model training. However, many mathematical programming problems in classical application areas are sparse, lack a uniform memory layout, and therefore do not benefit from the same kind of parallelism as dense linear algebra. As a result, integrating GPUs into mathematical programming solvers poses greater challenges and often necessitates substantial modifications to the overall algorithm.

On the one hand, first-order algorithms have emerged as a suitable class for GPU implementation. Since these algorithms rely on sparse matrix-vector multiplication and simple vector operations, implementing GPU acceleration is usually straightforward. Recent successful implementations include cuPDLP \cite{luCuPDLPCStrengthenedImplementation2024,lu2025cupdlpx}, cuOSQP \cite{schubigerGPUAccelerationADMM2020}, and cuOPT \cite{NVIDIACuopt2025}. However, the linear convergence rate of first-order methods restricts their effectiveness in applications requiring fast convergence, prompting the exploration of second-order alternatives for applications that require higher accuracy.

On the other hand, second-order solvers inherently rely on \emph{direct linear solvers}. For example, within \gls*{ipm}, each barrier iteration necessitates solving a linear system known as the \gls*{kkt} system. These systems become increasingly ill-conditioned as the iterate approaches the solution, rendering the use of iterative linear solvers, such as preconditioned Krylov methods, ineffective in most cases. Therefore, a reliable direct linear solver is a prerequisite for the effectiveness of second-order solvers. For years, the development of GPU-accelerated second-order solvers has been hindered by the absence of robust and efficient sparse direct linear solvers.

This status quo has changed with NVIDIA's release of cuDSS, a library of direct sparse linear solvers for GPUs \cite{nvidiaNVIDIACuDSSPreview}. It provides sparse Cholesky, LDL$^\top$, and LU factorization routines. While it currently lacks the LBL$^\top$ factorization capabilities commonly used for \gls*{nlp} solvers, its LDL$^\top$ and Cholesky functionalities are sufficient for implementing modified versions of the \gls*{ipm}. Consequently, cuDSS has spurred advances in GPU-accelerated second-order solvers, including MadNLP \cite{shinAcceleratingOptimalPower2024} and Clarabel \cite{goulartClarabelInteriorpointSolver2024}, achieving significant speedups on large-scale instances \cite{shinNVIDIACuDSSLibrary2024,shinAcceleratingOptimalPower2024,pacaudCondensedspaceMethodsNonlinear2024,shinScalableMultiPeriodAC2024,pacaudGPUacceleratedDynamicNonlinear2024}.

This paper provides an overview of the current state of the art in GPU implementations of second-order optimization solvers, with an emphasis on the following aspects: (i) The \gls*{ipm} is considered the primary mechanism for handling inequality constraints, as active-set methods are generally regarded as less scalable \cite{nocedalNumericalOptimization2006}. (ii) We mainly focus on solving KKT systems, since other components, such as line search and barrier updates, can be ported to GPUs straightforwardly using \texttt{map} or \texttt{reduce} operations. (iii) We primarily consider NVIDIA GPUs and the CUDA software stack, as they currently offer the most mature direct sparse solver implementation. (iv) Due to space constraints, hybrid \gls*{kkt} strategies \cite{regevHyKKTHybridDirectiterative2023}, reduced-space methods \cite{pacaudAcceleratingCondensedInteriorPoint2023}, and other domain-specific approaches \cite{adabagMPCGPURealTimeNonlinear2024} are not covered.

\section{Direct Linear Solvers for Optimization}\label{eqn:linear}
This section provides an overview of the direct linear algebra methods frequently employed in second-order methods and discusses the rationale behind the development of \emph{pivoting-free \gls*{ipm}}.

\paragraph{LDL$^\top$ Factorization.}
LDL$^\top$ factorization, a signed variant of Cholesky decomposition, decomposes a matrix $A$ into $LDL^\top$, where $L$ is lower triangular and $D$ is diagonal (for sparse systems, a fill-in reducing reordering $P$ must be employed, resulting in $P^\top A P = L D L^\top$). This method can be utilized to solve $Ax = b$, where the solution is obtained by first solving the lower triangular system $Ly = b$, followed by diagonal scaling with $D^{-1}$ and solving the upper triangular system $L^\top x = y$.

A notable property of LDL$^\top$ factorization is that, provided the matrix $A$ is \gls*{sqd}, the LDL$^\top$ factorization exists for any given permutation of the matrix (so-called \emph{strongly factorizable}) \cite{vanderbeiSymmetricQuasidefiniteMatrices1995}. This \emph{does not imply that numerical stability is guaranteed} for any reordering (see \cite{vanderbeiSymmetricQuasidefiniteMatrices1995}), but in practice, strong factorizability is often sufficient to ensure that these methods can be effectively utilized within optimization solvers \cite{stellatoOSQPOperatorSplitting2020}. Many \gls*{kkt} systems in optimization are \gls*{sqd}, can become \gls*{sqd} with infinitesimal regularization, or can be converted to \gls*{sqd} systems. If the system is \gls*{spd}, which is a sufficient condition for \gls*{sqd}, the LDL$^\top$ factorization or Cholesky factorization exists in a fill-in reducing manner, and the factorization process is always numerically stable. \gls*{spd} systems arise from unconstrained optimization problems or are obtained as a result of condensation, which will be discussed in \Cref{sec:ipm}.

\paragraph{Numerical Pivoting.}

For general indefinite matrices without \gls*{sqd} structure (e.g., augmented systems arising from nonconvex \glspl*{nlp}~\cite{wachterImplementationInteriorpointFilter2006}), the LDL$^\top$ factorization is not guaranteed to exist, and dynamic numerical pivoting is commonly employed to avoid zero pivots and improve the numerical stability of the factorization process. Dynamic numerical pivoting procedures examine a limited set of candidate pivots—typically within a row and column—and select the most suitable pivot according to a stability criterion \cite{schenkFASTFACTORIZATIONPIVOTING}. Three widely used dynamic pivoting strategies are Bunch–Kaufman, rook, and delayed pivoting, which select $1 \times 1$ or $2 \times 2$ pivots, although other variants and hybrid approaches also exist \cite{duffDirectMethodsSparse2017}. The variant of LDL$^\top$ with $2 \times 2$ pivots is often referred to as LBL$^\top$ factorization. If none of these methods succeed, the pivot is perturbed by a small value to allow numerical division \cite{schenkFASTFACTORIZATIONPIVOTING}. This procedure introduces numerical error, which must be corrected through iterative refinements. One of the drawbacks of numerical pivoting is that it requires deviating from the fill-in reducing reordering $P$, leading to additional fill-in and disrupting potential parallelism.

\paragraph{GPU Direct Solvers for Optimization}
As described above, the numerical pivoting procedure is crucial for ensuring the numerical stability of direct sparse linear solvers. However, implementing numerical pivoting has been recognized as one of the most challenging components of direct sparse linear solvers on GPUs, as these strategies are serial in nature \cite{swirydowiczLinearSolversPower2022}. Moreover, since coarse-grained tree-level parallelism must be employed to exploit GPU parallelism, numerical pivoting should be applied in a manner that does not disrupt the parallelism at the elimination tree level, further complicating the implementation. The current version of cuDSS has partial pivoting capabilities, but it does not support the LBL$^\top$ factorization implemented in CPU solvers \cite{nvidiaNVIDIACuDSSPreview}.

Therefore, to fully exploit the benefits of existing GPU direct solvers, it is crucial to ensure that \emph{the \gls*{kkt} system can be solved without numerical pivoting}, which motivates the development of \emph{pivoting-free interior-point methods}. This can be achieved by converting the \gls*{kkt} systems into an \gls*{sqd}, or even \gls*{spd} form, where strong factorizability guarantees the existence of the LDL$^\top$ factorization for any fill-in reducing reordering, allowing the factorization to succeed without relying on pivoting. This can be achieved through regularization or condensation, which we elaborate in \Cref{sec:ipm}. Once the pivoting requirement is eliminated, numerical factorization and triangular solves can be efficiently performed on GPUs \cite{naumovParallelSolutionSparse}. Although algorithms for computing fill-in reducing reorderings (e.g., minimum degree ordering \cite{amestoyApproximateMinimumDegree1996} or nested dissection \cite{karypisMETISSoftwarePackage1997}) are serial (e.g., cuDSS performs this operation on the CPU \cite{nvidiaNVIDIACuDSSPreview}), the reordering needs to be computed only once and can be reused, allowing the overhead to be amortized.

\section{Pivoting-Free Interior-Point Methods}\label{sec:ipm}
We now explain how the \gls*{ipm} can be adapted to avoid numerical pivoting, thereby enabling the use of GPU direct solvers relying only on static pivoting. We first provide a brief overview of the \gls*{ipm} and its KKT system formulation, followed by a discussion about condensed KKT systems.

\paragraph{Interior-Point Methods and KKT Systems.}
The \gls*{ipm} is a class of optimization algorithms designed to solve inequality-constrained optimization problems \cite{nocedalNumericalOptimization2006}.
In this work, we restrict attention to problems with inequality constraints only. This does not entail any loss of generality: equality constraints would introduce only minor and well-understood modifications to the \gls*{kkt} conditions and to the resulting linear systems.
The \gls*{ipm} transforms \cref{eqn:opt} into a sequence of log-barrier subproblems and attempts to solve its \gls*{kkt} conditions:
\begin{align}\label{eqn:kkt}
  \nabla f(x) - \nabla g(x)^\top \lambda = 0, \quad
  S \Lambda e - \mu e = 0, \quad
  g(x) - s = 0,
\end{align}
where $s \in \mathbb{R}^m$ denotes the slack variable used to reformulate the inequality constraints as equality constraints, $\mu > 0$ is the barrier parameter, $\lambda \in \mathbb{R}^m$ are the Lagrange multipliers, $S = \text{diag}(s)$, $\Lambda = \text{diag}(\lambda)$, and $e$ is the vector of ones.

The system~\cref{eqn:kkt} is solved using Newton's method. At each iteration, we obtain the search direction by solving the following (regularized and symmetric) KKT system:
\begin{align}\label{eqn:kkt_system}
  \begin{bmatrix}
    \nabla^2_{x x} \mathcal{L}(x,s,\lambda) + \delta_p I & & \nabla g(x)^\top
    \\ & S^{-1}\Lambda &  -I
    \\ \nabla g(x) & -I &  - \delta_d I
  \end{bmatrix}
  \begin{bmatrix}
    \phantom{-}d_x \\
    \phantom{-}d_s \\
    -d_\lambda
  \end{bmatrix} =
  -\begin{bmatrix}
    \nabla f(x) - \nabla g(x)^\top \lambda\\
    \Lambda e - \mu S^{-1}e \\
    g(x) - s\\
  \end{bmatrix},
\end{align}
where $\mathcal{L}(x,s,\lambda) := f(x) - \lambda^\top (g(x)-s)$ and $(\delta_p, \delta_d)$ are primal-dual regularization parameters.

\paragraph{Regularization.}
The regularization parameters $(\delta_p, \delta_d)$ ensure (i) the non-singularity of \cref{eqn:kkt_system} and/or (ii) the descent property of the Newton step. For convex problems, infinitesimal $\delta_p, \delta_d > 0$ ensures the \gls*{sqd} condition for the matrix in \cref{eqn:kkt_system}, implying strong factorizability. This idea has led to several robust \gls*{ipm} implementations on CPUs~\cite{friedlanderPrimalDualRegularized2012}. In nonconvex cases, primal-dual regularization provides a mechanism to get a descent direction for a given merit function \cite{wachterImplementationInteriorpointFilter2006}. \Gls*{ipm} solvers typically utilize a procedure known as \emph{inertia correction}, where the regularization parameters $(\delta_p, \delta_d)$ are increased until the number of positive, negative, and zero eigenvalues of \eqref{eqn:kkt_system} (collectively referred to as inertia, and available as a byproduct of the LDL$^\top$ and LBL$^\top$ factorizations) equals $(n+m, m, 0)$. Excessive regularization is undesirable, as it can potentially distort the step direction, leading to slow convergence.

\paragraph{Condensed KKT Systems.}
While \cref{eqn:kkt_system} is directly addressed by some solvers (e.g., Ipopt \cite{wachterImplementationInteriorpointFilter2006}), the system can be further \emph{condensed} into a so-called \emph{condensed \gls*{kkt} system}. In the context of GPU implementation, condensation offers advantages by either (i) reducing the system size and increasing its density—thereby providing more opportunities for parallelism—or (ii) enforcing the \gls*{spd} structure, which enables a pivoting-free implementation. However, depending on the sparsity pattern, the condensed system can become significantly denser, leading to higher memory requirements and computational overhead. Moreover, since the eliminated blocks are often highly ill-conditioned near the solution, the resulting condensed system may also suffer from ill-conditioning. Below, we outline several condensation strategies.
\begin{itemize}[leftmargin=*,itemsep=0pt,parsep=0pt,partopsep=0pt]
\item \textit{Augmented System}:
Since the \(S^{-1}\Lambda\) block in \cref{eqn:kkt_system} is always invertible due to the nature of the \gls*{ipm}, we can eliminate it to obtain the so-called \emph{augmented KKT system} \cite{nocedalNumericalOptimization2006}:
\begin{equation}\label{eqn:augmentedKKT}
  \begin{bmatrix}
    \nabla^2_{xx} \mathcal{L}(x,s,\lambda) + \delta_p I & \nabla g(x)^\top \\
    \nabla g(x) &  - \delta_d I - \Lambda^{-1} S
  \end{bmatrix}
  \begin{bmatrix}
    \phantom{-}d_x\\
    - d_\lambda
  \end{bmatrix} =
  -\begin{bmatrix}
    \nabla f(x) - \nabla g(x)^\top \lambda\\
    g(x) - \mu \Lambda^{-1} e
  \end{bmatrix}.
\end{equation}
This elimination does not incur significant computational overhead, and the number of non-zero entries in the resulting system does not increase.

\item \textit{Primal Condensed System}:
The \(\delta_d I + \Lambda^{-1}S\) block within \cref{eqn:augmentedKKT} is always invertible, and its elimination gives rise to a \emph{primal condensed KKT system}:
\begin{align}\label{eqn:kkt_primal}
  \left(\nabla^2_{xx} \mathcal{L}(x,s,\lambda) + \delta_p I + \nabla g(x)^\top (\delta_d I + \Lambda^{-1} S)^{-1} \nabla g(x) \right) d_x = - r_p \; ,
\end{align}
where \(r_p\) is an appropriate right-hand side derived from \cref{eqn:augmentedKKT}. This condensation has one key advantage for \glspl*{nlp}: the system becomes \gls*{spd} under the application of primal-dual regularization \((\delta_p, \delta_d)\) chosen based on the standard inertia correction procedure \cite{shinAcceleratingOptimalPower2024}, meaning that \emph{the system can be factorized using Cholesky factorization without numerical pivoting (or any reordering) in a numerically stable manner}. However, since the Jacobian \(\nabla g(x)\) can have dense rows, the condensed system can become arbitrarily dense, necessitating specialized treatment.

\item \textit{Dual Condensed System}:
When the problem is strongly convex or when the regularization parameter \(\delta_p\) is sufficiently large, the \(\nabla^2 \mathcal{L}(x,s,\lambda) + \delta_p I\) block is invertible, and by eliminating it, we obtain the \emph{dual condensed KKT system}:
\begin{align}\label{eqn:kkt_dual}
  \left(\delta_d I + \Lambda^{-1}S + \nabla g(x)\left(\nabla_{xx}^2 \mathcal{L}(x,s,\lambda) + \delta_p I\right)^{-1} \nabla g(x)^\top\right)
  d_\lambda = - r_d \; ,
\end{align}
where \(r_d\) is an appropriate right-hand side. The formulation in \eqref{eqn:kkt_dual} is often used as the default option for \gls*{lp} solvers with \(\delta_p>0\), and this system is often referred to as the \emph{normal equations}. Assuming that the primal Hessian is \gls*{spd}, this system is also \gls*{spd}, meaning that it can be stably factorized using Cholesky factorization without numerical pivoting. However, this system can also become arbitrarily dense when there is a dense column in \(\nabla g(x)\), which requires special treatment.
\end{itemize}

\paragraph{Pivoting-Free IPM.}
We now explain which KKT system formulation among \cref{eqn:kkt_system,eqn:augmentedKKT,eqn:kkt_dual,eqn:kkt_primal} is suitable for pivoting-free \gls*{ipm} implementations. The key requirement is that the KKT system matrix must be at least \gls*{sqd} without aggressive regularization. We detail the conditions below.
\begin{itemize}[leftmargin=*,itemsep=0pt,parsep=0pt,partopsep=0pt]
\item \textit{Convex Case}:
  For convex programs, all four formulations \cref{eqn:kkt_system,eqn:augmentedKKT,eqn:kkt_dual,eqn:kkt_primal} are appropriate, as any of these systems can become \gls*{sqd} for infinitesimal $\delta_p, \delta_d > 0$. However, \cref{eqn:kkt_dual,eqn:kkt_primal} may achieve better numerical stability due to their \gls*{spd} structure. MadIPM, an existing GPU \gls*{ipm} solver, employs \cref{eqn:kkt_system} with fixed primal-dual regularization.
\item \textit{Nonconvex Case}:
  For nonconvex problems, the primal condensed system \cref{eqn:kkt_primal} is the most suitable, as it can be made \gls*{spd} by choosing the primal-dual regularization parameters $(\delta_p, \delta_d)$ based solely on inertia correction. The augmented systems \cref{eqn:kkt_system,eqn:augmentedKKT} are not suitable because they are not guaranteed to be \gls*{sqd} unless aggressive (beyond what is necessary to ensure the descent condition) regularization parameters $(\delta_p, \delta_d)$ are used. The dual condensed system \cref{eqn:kkt_dual} is also unsuitable, as $(\nabla^2_{x x} \mathcal{L}(x,s,\lambda) + \delta_p I)^{-1}$ is difficult to compute due to nonlinear constraints. MadNLP, an existing GPU \gls*{nlp} solver, employs \cref{eqn:kkt_primal} with primal-dual regularization based on inertia correction.
\end{itemize}

\section{Algebraic Modeling Systems and Automatic Differentiation}\label{eqn:ad}
\Gls*{nlp} solvers require external oracles to evaluate $f$, $g$, and their first and second-order derivatives. In most modern optimization software stacks, the derivative evaluation code (either compiled or interpreted) is generated in a fully automated fashion through the so-called \emph{algebraic modeling systems}, which are typically equipped with \gls*{ad} capabilities, such as AMPL \cite{fourerModelingLanguageMathematical1990}, CasADi \cite{anderssonCasADiSoftwareFramework2019}, JuMP \cite{dunningJuMPModelingLanguage2017}, Pyomo \cite{hartPyomoModelingSolving2011}, and Gravity \cite{hijaziGravityMathematicalModeling2018}. As classical instances of mathematical programming problems are typically sparse, these systems have historically been developed independently of machine learning frameworks, which tend to focus more on dense problems.

To enable efficient derivative evaluations and ensure a fully GPU-resident optimization workflow, it is crucial to develop algebraic modeling systems that provide derivative evaluation code in the form of GPU kernels. To achieve this, one can concentrate on the observation that many practical instances of large-scale sparse mathematical programs exhibit highly repetitive structures. For example, $f$ may be a sum of many terms (e.g., $f(x) = \sum_{p\in P} \widetilde{f}(x; p)$), and $g$ may be a collection of numerous constraints generated from a common template (e.g., $g(x) = \left\{\widetilde{g}(x; p)\right\}_{p\in P}$). If such a structure exists, the evaluation and differentiation of $f$ and $g$ become embarrassingly parallel, making it feasible to construct GPU kernels for them.
This structure is particularly common in energy systems optimization, such as electricity or gas networks, and in engineering applications involving ODEs or PDEs, where constraints are derived from a discretization scheme (time or spatial).
In these cases, the same underlying equations are repeated over multiple nodes, time steps, or spatial elements, forming a natural template that GPU-parallel evaluation can exploit.
Emerging algebraic modeling systems, such as ExaModels.jl \cite{shinAcceleratingOptimalPower2024} or PyOptInterface \cite{yangPyOptInterfaceDesignImplementation2024}, are designed to capture this; for instance, ExaModels.jl requires users to specify the objective and constraint functions in the form of an iterator, such as
\begin{verbatim}
  objective(c, 100 * (x[i-1]^2 - x[i])^2 + (x[i-1] - 1)^2 for i = 2:N)
\end{verbatim}
which allows the user to inform the modeling system of repeated structures in the model. Then, the reverse-mode \gls*{ad} is applied to the template, and the resulting code is compiled into a GPU kernel. This approach enables efficient evaluation of the objective and constraints on GPUs, as well as the computation of their derivatives \cite{shinAcceleratingOptimalPower2024}.
It is worth emphasizing that the efficiency of GPU-based AD strongly depends on the presence of repeated structures in the problem.
Problems without repeated templates, where each constraint or objective term is unique, provide limited parallelism and thus only modest speed-ups on GPU.

\section{Numerical Results}\label{eqn:num}
We compared the performance of two GPU implementations (MadIPM for \glspl*{lp} and MadNLP for \glspl*{nlp}) against reference CPU solvers (Gurobi for \glspl*{lp} and Ipopt for \glspl*{nlp}).
\emph{We emphasize that the two solvers run entirely on the GPU, without data transfer between the host memory and the GPU memory.}
The results on the GPU have been generated using two Quadro GV100GPUs (with 128GB of memory). Despite being from an older generation, these two GPUs have very good double-precision performance (7.4 TFlops, compared to 14.8 TFlops for single precision performance).
We conducted the benchmark using MIPLIB 2010 (for \glspl*{lp}) \cite{kochMIPLIB20102011}, PGLIB-OPF (for \glspl*{nlp}) \cite{babaeinejadsarookolaeePowerGridLibrary2021}, and COPS (for \glspl*{nlp}) \cite{dolanBenchmarkingOptimizationSoftware2001}. The results are summarized in \Cref{tab:results}, and more details can be found in \Cref{apx:num}. These results can be reproduced using the source code available at \url{https://github.com/MadNLP/neurips2025-mathprog-on-gpu}. \textit{Disclaimer}: (i) The numerical results presented herein aim to demonstrate the current capabilities of GPU solvers by providing a comparison with comparable implementations on CPUs. This benchmark is not intended for a head-to-head performance comparison of the solvers. For example, some performance-critical options for CPU solvers, such as presolve and crossover, have been disabled to allow for a focused comparison of barrier iteration performance. Additionally, the convergence criteria for each solver differ slightly, and performance comparisons are based on user-facing tolerance options.

\begin{table}
  \footnotesize
  \begin{tabular}{|c|c|c|cc|cc|cc|cc|}
  \hline
  &\multirow{ 3}{*}{\bfseries Tol} & \multirow{ 3}{*}{\bfseries Solver} & \multicolumn{2}{c|}{\textbf{Small}}& \multicolumn{2}{c|}{\textbf{Medium}}& \multicolumn{2}{c|}{\textbf{Large}}& \multicolumn{2}{c|}{\multirow{2}{*}{\textbf{Total}}}\\
  &&& \multicolumn{2}{c|}{nnz $<2^{18}$}& \multicolumn{2}{c|}{$2^{18}\leq$ nnz $<2^{20}$}& \multicolumn{2}{c|}{$2^{20}\leq$ nnz}&&\\
  &&&  Solved & Time &  Solved & Time &  Solved & Time &  Solved & Time \\
  \hline\hline
  \multirow{4}{*}{\rotatebox{90}{\bfseries MIPLIB}}&    \multirow{2}{*}{$10^{-4}$} & MadIPM & 87 & 1.3013 & 56 & 5.0480 & 27 & 19.7925 & 170 & 4.5319  \\
  && Gurobi & 88 & 1.5439 & 58 & 10.4671 & 23 & 78.5783 & 169 & 9.3939  \\
  \cline{2-11}
  &\multirow{2}{*}{$10^{-8}$} & MadIPM & 85 & 2.8157 & 48 & 18.2642 & 25 & 33.1676 & 158 & 10.2820  \\
  && Gurobi & 88 & 1.5708 & 58 & 10.6148 & 24 & 76.3206 & 170 & 9.3826  \\
  \hline\hline
  \multirow{4}{*}{\rotatebox{90}{\bfseries OPF}}&\multirow{2}{*}{$10^{-4}$} & MadNLP & 31 & 0.4166 & 24 & 2.6380 & 11 & 3.7040 & 66 & 1.6979  \\
  && Ipopt & 31 & 0.3970 & 24 & 5.0697 & 11 & 38.5053 & 66 & 5.3817  \\
  \cline{2-11}
  &\multirow{2}{*}{$10^{-8}$} & MadNLP & 30 & 2.5037 & 24 & 4.6016 & 10 & 12.8040 & 64 & 4.6228  \\
  && Ipopt & 31 & 0.5100 & 24 & 5.4292 & 11 & 37.7818 & 66 & 5.5541  \\
  \hline\hline
  \multirow{4}{*}{\rotatebox{90}{\bfseries COPS}}&\multirow{2}{*}{$10^{-4}$} & MadNLP & 13 & 0.8665 & 15 & 4.8665 & 16 & 3.8194 & 44 & 3.2314  \\
  && Ipopt & 13 & 5.2315 & 15 & 15.9701 & 15 & 45.8411 & 43 & 19.2243  \\
  \cline{2-11}
  &\multirow{2}{*}{$10^{-8}$} & MadNLP & 13 & 0.8575 & 16 & 1.5572 & 16 & 8.3549 & 45 & 3.3797  \\
  && Ipopt & 13 & 5.9413 & 15 & 17.6758 & 15 & 40.8639 & 43 & 19.2999  \\
  \hline
\end{tabular}


  \caption{Solution times for CPU solvers (Gurobi and Ipopt) and GPU solvers (MadIPM and MadNLP) are represented using SGM10, defined as $(\prod_{i=1}^n (t_i + 10))^{1/n} - 10$, where $t_i$ denotes the solve time for the $i$-th instance (in seconds; unsolved instances are assigned a maximum wall time of 900 seconds) across various datasets: MIPLIB (88 small, 58 medium, and 28 large \glspl*{lp}), PGLIB-OPF (31 small, 24 medium, and 11 large \glspl*{nlp}), and COPS (13 small, 16 medium, and 16 large \glspl*{nlp}). For Gurobi, the Barrier method is used, with both the Presolve and Crossover options disabled. MadNLP is configured with cuDSS, while Ipopt is configured with either Ma27 (for PGLIB-OPF) or Ma57 (for COPS). All \glspl*{nlp} are modeled using ExaModels, which supports \gls*{nlp} function evaluation on both CPU and GPU. The benchmarking was conducted on a workstation equipped with two Intel Xeon Gold 6130 CPUs, two Quadro GV~100 GPUs, and 128 GB of memory.
  }
  \label{tab:results}
\end{table}

\paragraph{MIPLIB.}
We have performed the benchmark against a curated subset of instances within the MIPLIB 2010 library by selecting 174 instances that are sufficiently large and not trivially solved. We solve the linear programming relaxation associated to each problem (binary and integrality constraints are deleted). The results in \Cref{tab:results} indicate that the GPU solver can achieve, on average, approximately 4x speed-up for the 28 largest instances (with more than $2^{20}$ non-zeros) when the problems are solved to medium precision. The speed-up is relatively modest for medium-sized instances, and there is practically no advantage for small instances. This is expected, as the GPU solver is designed to handle large-scale problems, and small-scale problems cannot fully utilize the available parallel cores. In such cases, the overhead related to parallelism, such as task scheduling and thread launching, dominates the computation time rather than providing actual performance gains. For high precision, however, the speed-up is less pronounced, and the GPU solver solved significantly fewer instances.

\paragraph{PGLIB-OPF.}
We have benchmarked the performance of the solver for solving AC OPF problems based on polar power flow formulations \cite{PowerModelsJLOpenSource}. The results in \Cref{tab:results} indicate that the GPU solver can achieve an average speed-up of more than 10x for large instances when the problems are solved to medium precision. The speed-up is relatively modest for medium-sized instances, and there is practically no advantage for small instances. However, for high precision, the GPU solver does not reach the same level of robustness as the CPU solver, as the condensed system utilized by the GPU solvers often encounters worse conditioning; the GPU solver fails on two more instances. Nevertheless, the overall speed-up remains significant (3x on average for large instances).

\paragraph{COPS Benchmark.}
We conducted benchmarks using the COPS benchmark library on curated instances. The COPS benchmark instances are scalable, allowing users to specify the problem size. For each instance type, we formulated the problem in five different sizes, approximately doubling the number of variables and constraints each time. The results are similar to those of the PGLIB-OPF benchmark, but the speed-up is more pronounced in these instances. Again, for large instances, we can achieve more than a 10x speed-up on average at medium-precision ($10^{-4}$), but the speed-up
decreases to only 4x at high precision ($10^{-8}$).

\section{Conclusions and Future Outlook}\label{sec:conclusion}
We have presented an overview of the current landscape of GPU-accelerated second-order optimization solvers. With two specific existing solvers—MadIPM and MadNLP—and a modeling environment---ExaModels---we have demonstrated that GPU acceleration can achieve more than an order of magnitude speed-up for large instances when solved to medium precision. Solving problems robustly to high precision remains an open challenge, both for \gls*{lp} and \gls*{nlp} solvers.
Our upcoming GPU solver, MadNCL \cite{MadNCL}, dedicated to handling degenerate problems, aims to address the robustness issue we have observed at high precision.
Some open questions and implementation challenges are summarized below.

\begin{itemize}[leftmargin=*,itemsep=0pt,parsep=0pt,partopsep=0pt]
\item \textit{Numerical Precision of Condensed KKT Systems}: Condensed \gls*{kkt} systems are often preferred in pivoting-free implementations; however, stability can be compromised. Further research is needed to develop strategies to mitigate stability issues, especially for high-precision solves.
\item \textit{Batch Solvers}: GPU solvers can solve many small- and medium-sized problems in parallel, which can be facilitated through the implementation of batch solvers. Implementing second-order algorithms is feasible, as batch solutions (with or without uniform sparsity patterns) have been supported by CUDSS since version 0.6.
\item \textit{Hardware Portability}: Currently, most existing optimization and linear solvers are limited to NVIDIA GPUs. However, there is interest in developing hardware-agnostic solvers that can run on various GPU architectures, including AMD and Intel GPUs. A key requirement for this will be the development of portable sparse LDL$^\top$ factorizations.
\end{itemize}

\pagebreak

\bibliographystyle{plain}
\bibliography{shin}

\begin{thebibliography}{10}

\bibitem{PowerModelsJLOpenSource}
{{PowerModels}}. {{JL}}: {{An Open-Source Framework}} for {{Exploring Power
  Flow Formulations}} {\textbar} {{IEEE Conference Publication}} {\textbar}
  {{IEEE Xplore}}.
\newblock https://ieeexplore.ieee.org/abstract/document/8442948.

\bibitem{NVIDIACuopt2025}
{{NVIDIA}}/cuopt.
\newblock NVIDIA Corporation, August 2025.

\bibitem{adabagMPCGPURealTimeNonlinear2024}
Emre Adabag, Miloni Atal, William Gerard, and Brian Plancher.
\newblock {{MPCGPU}}: {{Real-Time Nonlinear Model Predictive Control}} through
  {{Preconditioned Conjugate Gradient}} on the {{GPU}}, March 2024.

\bibitem{amestoyApproximateMinimumDegree1996}
Patrick~R. Amestoy, Timothy~A. Davis, and Iain~S. Duff.
\newblock An {{Approximate Minimum Degree Ordering Algorithm}}.
\newblock {\em SIAM Journal on Matrix Analysis and Applications},
  17(4):886--905, October 1996.

\bibitem{anderssonCasADiSoftwareFramework2019}
Joel A.~E. Andersson, Joris Gillis, Greg Horn, James~B. Rawlings, and Moritz
  Diehl.
\newblock {{CasADi}}: A software framework for nonlinear optimization and
  optimal control.
\newblock {\em Mathematical Programming Computation}, 11(1):1--36, March 2019.

\bibitem{babaeinejadsarookolaeePowerGridLibrary2021}
Sogol Babaeinejadsarookolaee, Adam Birchfield, Richard~D. Christie, Carleton
  Coffrin, Christopher DeMarco, Ruisheng Diao, Michael Ferris, Stephane
  Fliscounakis, Scott Greene, Renke Huang, Cedric Josz, Roman Korab, Bernard
  Lesieutre, Jean Maeght, Terrence W.~K. Mak, Daniel~K. Molzahn, Thomas~J.
  Overbye, Patrick Panciatici, Byungkwon Park, Jonathan Snodgrass, Ahmad
  Tbaileh, Pascal~Van Hentenryck, and Ray Zimmerman.
\newblock The {{Power Grid Library}} for {{Benchmarking AC Optimal Power Flow
  Algorithms}}, January 2021.

\bibitem{dolanBenchmarkingOptimizationSoftware2001}
E.~D. Dolan and J.~J. More.
\newblock Benchmarking optimization software with {{COPS}}.
\newblock Technical Report ANL/MCS-TM-246, Argonne National Lab., IL (US),
  January 2001.

\bibitem{duffDirectMethodsSparse2017}
I.~S. Duff, A.~M. Erisman, and J.~K. Reid.
\newblock {\em Direct {{Methods}} for {{Sparse Matrices}}}.
\newblock Oxford University Press, March 2017.

\bibitem{dunningJuMPModelingLanguage2017}
Iain Dunning, Joey Huchette, and Miles Lubin.
\newblock {{JuMP}}: {{A Modeling Language}} for {{Mathematical Optimization}}.
\newblock {\em SIAM Review}, 59(2):295--320, January 2017.

\bibitem{fourerModelingLanguageMathematical1990}
Robert Fourer, David~M. Gay, and Brian~W. Kernighan.
\newblock A {{Modeling Language}} for {{Mathematical Programming}}.
\newblock {\em Management Science}, 36(5):519--554, May 1990.

\bibitem{friedlanderPrimalDualRegularized2012}
M.~P. Friedlander and D.~Orban.
\newblock A primal--dual regularized interior-point method for convex quadratic
  programs.
\newblock {\em Mathematical Programming Computation}, 4(1):71--107, March 2012.

\bibitem{goulartClarabelInteriorpointSolver2024}
Paul~J. Goulart and Yuwen Chen.
\newblock Clarabel: {{An}} interior-point solver for conic programs with
  quadratic objectives, May 2024.

\bibitem{hartPyomoModelingSolving2011}
William~E. Hart, Jean-Paul Watson, and David~L. Woodruff.
\newblock Pyomo: Modeling and solving mathematical programs in {{Python}}.
\newblock {\em Mathematical Programming Computation}, 3(3):219--260, September
  2011.

\bibitem{hijaziGravityMathematicalModeling2018}
Hassan Hijazi, Guanglei Wang, and Carleton Coffrin.
\newblock Gravity: {{A Mathematical Modeling Language}} for {{Optimization}}
  and {{Machine Learning}}.
\newblock October 2018.

\bibitem{karypisMETISSoftwarePackage1997}
George Karypis and Vipin Kumar.
\newblock {{METIS}}: {{A Software Package}} for {{Partitioning Unstructured
  Graphs}}, {{Partitioning Meshes}}, and {{Computing Fill-Reducing Orderings}}
  of {{Sparse Matrices}}.
\newblock 1997.

\bibitem{kochMIPLIB20102011}
Thorsten Koch, Tobias Achterberg, Erling Andersen, Oliver Bastert, Timo
  Berthold, Robert~E. Bixby, Emilie Danna, Gerald Gamrath, Ambros~M. Gleixner,
  Stefan Heinz, Andrea Lodi, Hans Mittelmann, Ted Ralphs, Domenico Salvagnin,
  Daniel~E. Steffy, and Kati Wolter.
\newblock {{MIPLIB}} 2010.
\newblock {\em Mathematical Programming Computation}, 3(2):103--163, June 2011.

\bibitem{lu2025cupdlpx}
Haihao Lu, Zedong Peng, and Jinwen Yang.
\newblock {cuPDLPx: A further enhanced GPU-based first-order solver for linear
  programming}.
\newblock {\em arXiv preprint arXiv:2507.14051}, 2025.

\bibitem{luCuPDLPCStrengthenedImplementation2024}
Haihao Lu, Jinwen Yang, Haodong Hu, Qi~Huangfu, Jinsong Liu, Tianhao Liu, Yinyu
  Ye, Chuwen Zhang, and Dongdong Ge.
\newblock {{cuPDLP-C}}: {{A Strengthened Implementation}} of {{cuPDLP}} for
  {{Linear Programming}} by {{C}} language, January 2024.

\bibitem{MadNCL}
Alexis Montoison, François Pacaud, Michael Saunders, Sungho Shin, and
  Dominique Orban.
\newblock {MadNCL: A GPU Implementation of Algorithm NCL for Large-Scale,
  Degenerate Nonlinear Programs}.
\newblock {\em arXiv preprint arXiv:2510.05885}, 2025.

\bibitem{naumovParallelSolutionSparse}
Maxim Naumov.
\newblock Parallel {{Solution}} of {{Sparse Triangular Linear Systems}} in the
  {{Preconditioned Iterative Methods}} on the {{GPU}}.

\bibitem{nocedalNumericalOptimization2006}
Jorge Nocedal and Stephen~J. Wright.
\newblock {\em Numerical Optimization}.
\newblock Springer Series in Operations Research. Springer, New York, 2nd ed
  edition, 2006.

\bibitem{nvidiaNVIDIACuDSSPreview}
NVIDIA.
\newblock {{NVIDIA cuDSS}} ({{Preview}}): {{A}} high-performance {{CUDA
  Library}} for {{Direct Sparse Solvers}} --- {{NVIDIA cuDSS}} documentation.
\newblock https://docs.nvidia.com/cuda/cudss/index.html.

\bibitem{pacaudGPUacceleratedDynamicNonlinear2024}
Fran{\c c}ois Pacaud and Sungho Shin.
\newblock {{GPU-accelerated}} dynamic nonlinear optimization with {{ExaModels}}
  and {{MadNLP}}.
\newblock In {\em 2024 {{IEEE}} 63rd {{Conference}} on {{Decision}} and
  {{Control}} ({{CDC}})}, pages 5963--5968, December 2024.

\bibitem{pacaudCondensedspaceMethodsNonlinear2024}
Fran{\c c}ois Pacaud, Sungho Shin, Alexis Montoison, Michel Schanen, and Mihai
  Anitescu.
\newblock Condensed-space methods for nonlinear programming on {{GPUs}}, May
  2024.

\bibitem{pacaudAcceleratingCondensedInteriorPoint2023}
Fran{\c c}ois Pacaud, Sungho Shin, Michel Schanen, Daniel~Adrian Maldonado, and
  Mihai Anitescu.
\newblock Accelerating {{Condensed Interior-Point Methods}} on {{SIMD}}/{{GPU
  Architectures}}.
\newblock {\em Journal of Optimization Theory and Applications}, February 2023.

\bibitem{regevHyKKTHybridDirectiterative2023}
Shaked Regev, Nai-Yuan Chiang, Eric Darve, Cosmin~G. Petra, Michael~A.
  Saunders, Kasia {\'S}wirydowicz, and Slaven Pele{\v s}.
\newblock {{HyKKT}}: A hybrid direct-iterative method for solving {{KKT}}
  linear systems.
\newblock {\em Optimization Methods and Software}, 38(2):332--355, March 2023.

\bibitem{schenkFASTFACTORIZATIONPIVOTING}
Olaf Schenk.
\newblock {{ON FAST FACTORIZATION PIVOTING METHODS FOR SPARSE SYMMETRIC
  INDEFINITE SYSTEMS}}.

\bibitem{schubigerGPUAccelerationADMM2020}
Michel Schubiger, Goran Banjac, and John Lygeros.
\newblock {{GPU}} acceleration of {{ADMM}} for large-scale quadratic
  programming.
\newblock {\em Journal of Parallel and Distributed Computing}, 144:55--67,
  October 2020.

\bibitem{shinAcceleratingOptimalPower2024}
Sungho Shin, Mihai Anitescu, and Fran{\c c}ois Pacaud.
\newblock Accelerating optimal power flow with {{GPUs}}: {{SIMD}} abstraction
  of nonlinear programs and condensed-space interior-point methods.
\newblock {\em Electric Power Systems Research}, 236:110651, November 2024.

\bibitem{shinNVIDIACuDSSLibrary2024}
Sungho Shin, Fran{\c c}ois Pacaud, Alexis Montoison, Mark Wolf, and Becca
  Zandstein.
\newblock {{NVIDIA cuDSS Library Removes Barriers}} to {{Optimizing}} the {{US
  Power Grid}}.
\newblock
  https://developer.nvidia.com/blog/nvidia-cudss-library-removes-barriers-to-optimizing-the-us-power-grid/,
  November 2024.

\bibitem{shinScalableMultiPeriodAC2024}
Sungho Shin, Vishwas Rao, Michel Schanen, D.~Adrian Maldonado, and Mihai
  Anitescu.
\newblock Scalable {{Multi-Period AC Optimal Power Flow Utilizing GPUs}} with
  {{High Memory Capacities}}, May 2024.

\bibitem{stellatoOSQPOperatorSplitting2020}
Bartolomeo Stellato, Goran Banjac, Paul Goulart, Alberto Bemporad, and Stephen
  Boyd.
\newblock {{OSQP}}: An operator splitting solver for quadratic programs.
\newblock {\em Mathematical Programming Computation}, 12(4):637--672, December
  2020.

\bibitem{swirydowiczLinearSolversPower2022}
Kasia {\'S}wirydowicz, Eric Darve, Wesley Jones, Jonathan Maack, Shaked Regev,
  Michael~A. Saunders, Stephen~J. Thomas, and Slaven Pele{\v s}.
\newblock Linear solvers for power grid optimization problems: {{A}} review of
  {{GPU-accelerated}} linear solvers.
\newblock {\em Parallel Computing}, 111:102870, July 2022.

\bibitem{vanderbeiSymmetricQuasidefiniteMatrices1995}
Robert~J. Vanderbei.
\newblock Symmetric {{Quasidefinite Matrices}}.
\newblock {\em SIAM Journal on Optimization}, 5(1):100--113, February 1995.

\bibitem{wachterImplementationInteriorpointFilter2006}
Andreas W{\"a}chter and Lorenz~T. Biegler.
\newblock On the implementation of an interior-point filter line-search
  algorithm for large-scale nonlinear programming.
\newblock {\em Mathematical Programming}, 106(1):25--57, March 2006.

\bibitem{yangPyOptInterfaceDesignImplementation2024}
Yue Yang, Chenhui Lin, Luo Xu, and Wenchuan Wu.
\newblock {{PyOptInterface}}: {{Design}} and implementation of an efficient
  modeling language for mathematical optimization, May 2024.

\end{thebibliography}

\appendix
\crefalias{section}{appendix}

\section{More Details on Numerical Results}\label{apx:num}
\subsection{Solver Options}
\subsection{Gurobi}
\begin{verbatim}
FeasibilityTol = 1e-4 or 1e-8
OptimalityTol = 1e-4 or 1e-8
TimeLimit = 900.0
Method = 2
Presolve = 0
Crossover = 0
Threads = 16
\end{verbatim}
\subsection{Ipopt}
\begin{verbatim}
tol = 1e-4 or 1e-8
bound_relax_factor = 1e-4 or 1e-8
max_wall_time = 900.0
linear_solver = "ma27" or "ma57"
ma57_automatic_scaling = "yes"
dual_inf_tol = 10000.0
constr_viol_tol = 10000.0
compl_inf_tol = 10000.0
honor_original_bounds = "no"
print_timing_statistics = "yes"
\end{verbatim}
\subsection{MadIPM}
\begin{verbatim}
tol = 1e-4 or 1e-8
max_wall_time = 900.0
max_iter = 500
linear_solver = MadNLPGPU.CUDSSSolver
cudss_algorithm = MadNLP.LDL
regularization = MadIPM.FixedRegularization(1e-8, -1e-8)
print_level = MadNLP.INFO
rethrow_error = true
\end{verbatim}
\subsection{MadNLP}
\begin{verbatim}
tol = 1e-4 or 1e-8
max_wall_time = 900.0
\end{verbatim}
\subsection{Full Numerical Results}
\begin{verbatim}
--------------------------------------------------------------------
              MIPLIB benchmark results (tol = 0.0001)
--------------------------------------------------------------------
           problem | log2(nnz)|      MadIPM      |      Gurobi        
                   |          | solved|     time | solved|     time     
--------------------------------------------------------------------
              n3-3 |    15.15 |   1   |     0.30 |   1   |     0.26
       neos-506422 |    15.26 |   1   |     0.14 |   1   |     0.15
            ramos3 |    15.26 |   1   |     0.30 |   1   |     0.39
      iis-bupa-cov |    15.42 |   1   |     0.22 |   1   |     0.27
       neos-777800 |    15.47 |   1   |     0.22 |   1   |     0.17
            d10200 |    15.54 |   1   |     0.24 |   1   |     0.16
            hanoi5 |    15.69 |   1   |     0.35 |   1   |     0.78
         ns1778858 |    15.69 |   1   |     0.22 |   1   |     0.28
           eil33-2 |    15.70 |   1   |     0.14 |   1   |     0.17
       neos-941262 |    15.76 |   1   |     0.40 |   1   |     0.42
   lectsched-4-obj |    15.77 |   1   |     0.19 |   1   |     0.24
       neos-984165 |    15.79 |   1   |     0.41 |   1   |     0.44
       neos-935769 |    15.80 |   1   |     0.32 |   1   |     0.37
       neos-948126 |    15.85 |   1   |     0.39 |   1   |     0.44
        reblock166 |    15.87 |   1   |     0.51 |   1   |     2.31
           lrsa120 |    15.91 |   1   |     0.21 |   1   |     0.79
       neos-935627 |    15.95 |   1   |     0.40 |   1   |     0.44
  rococoC12-111000 |    15.96 |   1   |     0.83 |   1   |     1.50
      neos-1171737 |    15.99 |   1   |     0.29 |   1   |     0.21
       neos-937511 |    16.08 |   1   |     0.43 |   1   |     0.44
            sp98ir |    16.14 |   1   |     0.25 |   1   |     0.22
         atm20-100 |    16.18 |   0   |     0.13 |   1   |     0.43
    methanosarcina |    16.18 |   1   |     0.21 |   1   |     0.95
       neos-937815 |    16.19 |   1   |     0.54 |   1   |     0.53
       neos-826812 |    16.22 |   1   |     0.57 |   1   |     0.38
    satellites1-25 |    16.26 |   1   |     0.86 |   1   |     1.14
          wachplan |    16.27 |   1   |     0.28 |   1   |     0.27
      iis-pima-cov |    16.30 |   1   |     0.35 |   1   |     0.82
          dano3mip |    16.34 |   1   |     0.85 |   1   |     4.34
           biella1 |    16.35 |   1   |     0.75 |   1   |     0.46
           30n20b8 |    16.39 |   1   |     0.37 |   1   |     0.44
             air04 |    16.47 |   1   |     0.41 |   1   |     0.47
       neos-826694 |    16.52 |   1   |     0.50 |   1   |     0.35
         queens-30 |    16.55 |   1   |     0.14 |   1   |     0.22
      neos-1605075 |    16.64 |   1   |     0.77 |   1   |     1.15
      neos-1605061 |    16.66 |   1   |     0.99 |   1   |     1.27
   ash608gpia-3col |    16.68 |   1   |     0.64 |   1   |     1.31
          blp-ic97 |    16.74 |   1   |     0.32 |   1   |     0.28
            sts405 |    16.75 |   1   |     0.36 |   1   |     0.39
             sct32 |    16.77 |   1   |     0.74 |   1   |     0.86
        opm2-z7-s2 |    16.82 |   1   |     2.11 |   1   |     5.64
      rmatr200-p20 |    16.85 |   1   |     0.50 |   1   |     1.86
       neos-693347 |    16.86 |   1   |     0.47 |   1   |     0.64
       lectsched-2 |    16.87 |   1   |     0.37 |   1   |     0.55
      neos-1109824 |    16.89 |   1   |     0.74 |   1   |    12.48
              sct1 |    16.89 |   1   |     1.70 |   1   |     1.68
             net12 |    16.90 |   1   |     1.00 |   1   |     3.85
         momentum1 |    17.04 |   1   |     0.99 |   1   |    14.44
         shipsched |    17.06 |   1   |     0.77 |   1   |     0.87
              dc1c |    17.07 |   1   |     1.16 |   1   |     1.12
       neos-916792 |    17.07 |   1   |     0.21 |   1   |     0.71
              leo1 |    17.08 |   1   |     0.51 |   1   |     0.90
      rmatr200-p10 |    17.10 |   1   |     0.55 |   1   |     2.35
       neos-738098 |    17.14 |   1   |     0.57 |   1   |     0.68
       rmatr200-p5 |    17.20 |   1   |     0.55 |   1   |     3.45
       neos-952987 |    17.22 |   1   |     0.54 |   1   |     0.81
         ex1010-pi |    17.23 |   1   |     0.99 |   1   |     1.08
            mzzv11 |    17.30 |   1   |     1.15 |   1   |     2.14
       neos-934278 |    17.43 |   1   |     1.13 |   1   |     1.16
        neos808444 |    17.44 |   1   |     0.84 |   1   |     0.82
            d20200 |    17.44 |   1   |     0.41 |   1   |     0.43
       lectsched-3 |    17.48 |   1   |     0.59 |   1   |     0.86
              sct5 |    17.50 |   1   |     1.34 |   1   |     2.19
          germanrr |    17.51 |   1   |     0.44 |   1   |     0.91
 satellites2-60-fs |    17.59 |   1   |     1.61 |   1   |     2.37
              bab5 |    17.63 |   1   |     0.88 |   1   |     0.81
       lectsched-1 |    17.64 |   1   |     0.59 |   1   |     0.96
   lectsched-1-obj |    17.64 |   1   |     0.61 |   1   |     1.13
       neos-824661 |    17.65 |   1   |     1.16 |   1   |     0.70
             t1722 |    17.65 |   1   |     0.66 |   1   |     0.89
      core2536-691 |    17.68 |   1   |     1.03 |   1   |     2.20
            dolom1 |    17.68 |   1   |     1.28 |   1   |     1.11
          ns930473 |    17.68 |   1   |     0.99 |   1   |     1.14
        reblock420 |    17.68 |   1   |     2.31 |   1   |    18.20
             sing2 |    17.70 |   1   |     0.95 |   1   |     1.95
         ns1456591 |    17.71 |   1   |     0.77 |   1   |     0.35
          blp-ar98 |    17.73 |   1   |     0.53 |   1   |     0.49
         stockholm |    17.81 |   1   |     4.06 |   1   |     9.09
              leo2 |    17.82 |   1   |     0.68 |   1   |     0.45
              bab1 |    17.87 |   1   |     0.84 |   1   |     1.23
           rmine10 |    17.90 |   1   |     2.94 |   1   |    28.49
       neos-933966 |    17.92 |   1   |     1.83 |   1   |     1.19
         uc-case11 |    17.92 |   1   |     1.69 |   1   |     2.26
        rocII-4-11 |    17.96 |   1   |     0.76 |   1   |     0.83
       neos-933638 |    17.98 |   1   |     1.83 |   1   |     1.42
       neos-885086 |    18.01 |   1   |     0.67 |   1   |     0.61
            app1-2 |    18.01 |   1   |     0.92 |   1   |     4.54
             neos6 |    18.04 |   1   |     0.59 |   1   |     0.64
     core4872-1529 |    18.05 |   1   |     2.16 |   1   |     4.09
         ns1685374 |    18.07 |   1   |     1.38 |   1   |     2.74
            neos13 |    18.09 |   1   |     0.55 |   1   |     0.47
            sp97ar |    18.29 |   1   |     0.74 |   1   |     0.63
         ns2124243 |    18.30 |   1   |     1.34 |   1   |     1.16
         ns1905797 |    18.32 |   1   |     1.35 |   1   |    51.37
         ns1952667 |    18.36 |   1   |     0.21 |   1   |     0.13
            sp98ic |    18.37 |   1   |     0.59 |   1   |     0.59
       tanglegram1 |    18.39 |   1   |     0.83 |   1   |     1.04
         momentum2 |    18.43 |   1   |     1.14 |   1   |     7.72
        atlanta-ip |    18.44 |   0   |    14.87 |   1   |    17.60
          circ10-3 |    18.44 |   1   |     0.97 |   1   |    14.49
            sts729 |    18.44 |   1   |     2.48 |   1   |     1.16
             map14 |    18.45 |   1   |     4.22 |   1   |    11.68
             map20 |    18.45 |   1   |     4.36 |   1   |     9.01
             map06 |    18.45 |   1   |     4.19 |   1   |    12.88
             map10 |    18.45 |   1   |     4.14 |   1   |    11.30
          uc-case3 |    18.45 |   1   |     1.69 |   1   |     3.44
             map18 |    18.45 |   1   |     4.58 |   1   |     8.98
    satellites2-60 |    18.46 |   1   |     4.30 |   1   |     6.45
       neos-520729 |    18.47 |   1   |     1.25 |   1   |     2.15
       neos-957389 |    18.51 |   1   |     0.77 |   1   |     1.08
             nsr8k |    18.73 |   0   |     5.51 |   1   |     5.43
       neos-885524 |    18.75 |   1   |     1.37 |   1   |     1.19
 satellites3-40-fs |    18.80 |   1   |     4.34 |   1   |    10.44
        rocII-7-11 |    18.83 |   1   |     1.36 |   1   |     1.39
           vpphard |    18.85 |   1   |     2.07 |   1   |     6.26
             t1717 |    18.85 |   1   |     1.87 |   1   |     2.22
               ex9 |    18.93 |   1   |     1.75 |   1   |    22.39
               van |    18.99 |   1   |     1.80 |   1   |     5.37
              dc1l |    18.99 |   1   |     1.63 |   1   |     2.63
       opm2-z10-s2 |    19.05 |   1   |     6.91 |   1   |   387.96
          triptim1 |    19.08 |   1   |     2.80 |   1   |     7.54
          triptim2 |    19.09 |   1   |     3.41 |   1   |     7.50
       neos-506428 |    19.09 |   1   |     1.81 |   1   |     4.44
       neos-932816 |    19.09 |   1   |     2.21 |   1   |     2.81
          triptim3 |    19.10 |   1   |     3.24 |   1   |     6.72
         ns1116954 |    19.11 |   1   |    20.40 |   1   |   468.65
         ns1904248 |    19.18 |   1   |     1.90 |   1   |     9.19
           rail507 |    19.18 |   1   |     2.92 |   1   |     2.45
         ns1111636 |    19.19 |   1   |     1.92 |   1   |     1.82
        rocII-9-11 |    19.20 |   1   |     1.87 |   1   |     1.85
             n15-3 |    19.23 |   1   |     2.83 |   1   |     2.90
            rail01 |    19.26 |   1   |     7.10 |   1   |     9.42
        gmut-75-50 |    19.30 |   1   |     2.28 |   1   |     2.52
   pb-simp-nonunif |    19.41 |   1   |     1.93 |   1   |     4.16
       opm2-z11-s8 |    19.52 |   1   |    10.95 |   1   |   395.59
       neos-941313 |    19.67 |   1   |     3.87 |   1   |     2.76
    satellites3-40 |    19.73 |   1   |    17.34 |   1   |    18.92
       neos-859770 |    19.76 |   1   |     0.83 |   1   |     1.03
      neos-1140050 |    19.76 |   1   |     1.15 |   1   |    62.57
      netdiversion |    19.89 |   1   |     5.16 |   1   |     9.36
           rmine14 |    19.92 |   1   |    17.69 |   1   |   177.85
         momentum3 |    19.98 |   1   |     3.14 |   1   |   136.43
    buildingenergy |    19.99 |   1   |     6.55 |   1   |    21.72
           rvb-sub |    20.00 |   1   |     1.81 |   1   |     3.16
      opm2-z12-s14 |    20.02 |   1   |    13.01 |   0   |   905.66
       opm2-z12-s7 |    20.02 |   1   |    16.17 |   0   |   903.53
          vpphard2 |    20.03 |   1   |     5.81 |   1   |    23.73
             stp3d |    20.05 |   1   |    13.54 |   1   |    22.59
         eilA101-2 |    20.06 |   1   |     2.27 |   1   |     3.56
            npmv07 |    20.07 |   0   |    47.94 |   1   |    12.03
         ns2118727 |    20.10 |   1   |     7.03 |   1   |    31.67
           sing245 |    20.10 |   1   |     8.00 |   1   |   100.44
         ns2137859 |    20.11 |   1   |     7.51 |   1   |     5.41
              ex10 |    20.17 |   1   |     5.01 |   1   |   136.32
         ns1854840 |    20.27 |   1   |     6.48 |   1   |    13.41
            rail02 |    20.31 |   1   |    18.35 |   1   |    32.93
       neos-631710 |    20.35 |   1   |     4.41 |   1   |     6.49
           datt256 |    20.53 |   1   |     8.10 |   0   |   900.34
         ns1758913 |    20.89 |   1   |    11.03 |   1   |   886.11
      neos-1429212 |    20.93 |   1   |     4.19 |   1   |    11.67
  wnq-n100-mw99-14 |    20.94 |   1   |    25.02 |   1   |   809.00
         ns1853823 |    20.94 |   1   |    13.41 |   1   |    82.23
            co-100 |    21.00 |   1   |     3.00 |   1   |     4.28
            rail03 |    21.63 |   1   |    48.47 |   1   |    98.37
           n3seq24 |    21.73 |   1   |     7.47 |   1   |    15.93
       neos-476283 |    21.92 |   1   |    46.87 |   1   |    29.43
              bab3 |    21.97 |   1   |    22.44 |   1   |    24.84
           rmine21 |    22.33 |   1   |   233.02 |   0   |   922.05
         ivu06-big |    24.72 |   1   |   337.38 |   1   |   150.96
            mspp16 |    24.74 |   1   |    57.65 |   1   |   544.68
--------------------------------------------------------------------
\end{verbatim}

\begin{verbatim}
--------------------------------------------------------------------
              MIPLIB benchmark results (tol = 1.0e-8)
--------------------------------------------------------------------
           problem | log2(nnz)|      MadIPM      |      Gurobi        
                   |          | solved|     time | solved|     time     
--------------------------------------------------------------------
              n3-3 |    15.15 |   1   |     0.31 |   1   |     0.23
       neos-506422 |    15.26 |   1   |     0.20 |   1   |     0.10
            ramos3 |    15.26 |   1   |     0.32 |   1   |     0.42
      iis-bupa-cov |    15.42 |   1   |     0.22 |   1   |     0.28
       neos-777800 |    15.47 |   1   |     0.22 |   1   |     0.15
            d10200 |    15.54 |   1   |     0.28 |   1   |     0.17
         ns1778858 |    15.69 |   1   |     3.59 |   1   |     0.31
            hanoi5 |    15.69 |   1   |     0.33 |   1   |     0.76
           eil33-2 |    15.70 |   1   |     0.15 |   1   |     0.14
       neos-941262 |    15.76 |   1   |     0.44 |   1   |     0.41
   lectsched-4-obj |    15.77 |   1   |     0.21 |   1   |     0.23
       neos-984165 |    15.79 |   1   |     0.45 |   1   |     0.43
       neos-935769 |    15.80 |   1   |     0.37 |   1   |     0.37
       neos-948126 |    15.85 |   1   |     0.38 |   1   |     0.42
        reblock166 |    15.87 |   1   |     0.94 |   1   |     2.40
           lrsa120 |    15.91 |   1   |     0.26 |   1   |     0.36
       neos-935627 |    15.95 |   1   |     0.44 |   1   |     0.45
  rococoC12-111000 |    15.96 |   1   |     0.98 |   1   |     1.56
      neos-1171737 |    15.99 |   1   |     0.30 |   1   |     0.21
       neos-937511 |    16.08 |   1   |     0.42 |   1   |     0.41
            sp98ir |    16.14 |   1   |     0.27 |   1   |     0.20
    methanosarcina |    16.18 |   1   |     0.33 |   1   |     5.63
         atm20-100 |    16.18 |   0   |     0.14 |   1   |     0.41
       neos-937815 |    16.19 |   1   |     0.52 |   1   |     0.53
       neos-826812 |    16.22 |   1   |     0.61 |   1   |     0.39
    satellites1-25 |    16.26 |   1   |     7.92 |   1   |     1.06
          wachplan |    16.27 |   1   |     0.31 |   1   |     0.26
      iis-pima-cov |    16.30 |   1   |     0.37 |   1   |     0.76
          dano3mip |    16.34 |   1   |     1.06 |   1   |     1.48
           biella1 |    16.35 |   1   |     0.79 |   1   |     0.46
           30n20b8 |    16.39 |   1   |     0.41 |   1   |     0.45
             air04 |    16.47 |   1   |     0.50 |   1   |     0.40
       neos-826694 |    16.52 |   1   |     0.51 |   1   |     0.34
         queens-30 |    16.55 |   1   |     0.14 |   1   |     0.23
      neos-1605075 |    16.64 |   1   |     1.00 |   1   |     1.22
      neos-1605061 |    16.66 |   1   |     0.96 |   1   |     1.28
   ash608gpia-3col |    16.68 |   1   |     0.72 |   1   |     2.15
          blp-ic97 |    16.74 |   1   |     0.30 |   1   |     0.29
            sts405 |    16.75 |   1   |     0.38 |   1   |     0.37
             sct32 |    16.77 |   0   |     5.31 |   1   |     0.86
        opm2-z7-s2 |    16.82 |   1   |     2.20 |   1   |     5.70
      rmatr200-p20 |    16.85 |   1   |     0.93 |   1   |     1.97
       neos-693347 |    16.86 |   1   |     0.49 |   1   |     0.59
       lectsched-2 |    16.87 |   1   |     0.41 |   1   |     0.57
              sct1 |    16.89 |   0   |    10.24 |   1   |     1.65
      neos-1109824 |    16.89 |   1   |     0.73 |   1   |    12.76
             net12 |    16.90 |   1   |     1.01 |   1   |     3.76
         momentum1 |    17.04 |   1   |     1.44 |   1   |    14.47
         shipsched |    17.06 |   1   |     0.60 |   1   |     0.88
       neos-916792 |    17.07 |   1   |     0.22 |   1   |     0.27
              dc1c |    17.07 |   1   |     4.77 |   1   |     1.13
              leo1 |    17.08 |   1   |     0.37 |   1   |     0.26
      rmatr200-p10 |    17.10 |   1   |     1.08 |   1   |     2.38
       neos-738098 |    17.14 |   1   |     0.69 |   1   |     0.72
       rmatr200-p5 |    17.20 |   1   |     1.13 |   1   |     3.07
       neos-952987 |    17.22 |   1   |     0.67 |   1   |     0.85
         ex1010-pi |    17.23 |   1   |     1.21 |   1   |     1.14
            mzzv11 |    17.30 |   1   |     1.41 |   1   |     1.89
       neos-934278 |    17.43 |   1   |     1.23 |   1   |     1.11
        neos808444 |    17.44 |   1   |     0.82 |   1   |     0.76
            d20200 |    17.44 |   1   |     0.47 |   1   |     0.46
       lectsched-3 |    17.48 |   1   |     0.53 |   1   |     0.88
              sct5 |    17.50 |   1   |     4.20 |   1   |     2.15
          germanrr |    17.51 |   1   |     0.48 |   1   |     0.93
 satellites2-60-fs |    17.59 |   1   |     4.30 |   1   |     2.60
              bab5 |    17.63 |   1   |     1.10 |   1   |     0.81
       lectsched-1 |    17.64 |   1   |     0.55 |   1   |     0.94
   lectsched-1-obj |    17.64 |   1   |     0.64 |   1   |     1.24
             t1722 |    17.65 |   1   |     0.78 |   1   |     0.89
       neos-824661 |    17.65 |   1   |     1.26 |   1   |     0.75
      core2536-691 |    17.68 |   1   |     1.34 |   1   |     2.27
          ns930473 |    17.68 |   1   |     1.08 |   1   |     1.22
        reblock420 |    17.68 |   1   |     2.17 |   1   |    18.07
            dolom1 |    17.68 |   1   |     1.33 |   1   |     1.13
             sing2 |    17.70 |   1   |     1.24 |   1   |     1.96
         ns1456591 |    17.71 |   1   |     0.94 |   1   |     0.39
          blp-ar98 |    17.73 |   1   |     0.60 |   1   |     0.49
         stockholm |    17.81 |   1   |     3.82 |   1   |     9.31
              leo2 |    17.82 |   1   |     1.05 |   1   |     0.94
              bab1 |    17.87 |   1   |     0.90 |   1   |     1.21
           rmine10 |    17.90 |   1   |     3.33 |   1   |    28.39
         uc-case11 |    17.92 |   1   |     1.89 |   1   |     2.28
       neos-933966 |    17.92 |   1   |     1.96 |   1   |     1.21
        rocII-4-11 |    17.96 |   1   |     0.90 |   1   |     0.86
       neos-933638 |    17.98 |   1   |     1.83 |   1   |     1.42
       neos-885086 |    18.01 |   1   |     0.71 |   1   |     0.61
            app1-2 |    18.01 |   1   |     1.19 |   1   |     4.37
             neos6 |    18.04 |   1   |     0.59 |   1   |     0.61
     core4872-1529 |    18.05 |   1   |     2.95 |   1   |     4.87
         ns1685374 |    18.07 |   1   |     1.44 |   1   |     3.81
            neos13 |    18.09 |   1   |     0.67 |   1   |     0.52
            sp97ar |    18.29 |   1   |     0.93 |   1   |     0.64
         ns2124243 |    18.30 |   1   |     1.49 |   1   |     1.21
         ns1905797 |    18.32 |   1   |     1.44 |   1   |    49.25
         ns1952667 |    18.36 |   1   |     0.52 |   1   |     0.13
            sp98ic |    18.37 |   1   |     0.63 |   1   |     0.68
       tanglegram1 |    18.39 |   1   |     0.91 |   1   |     0.99
         momentum2 |    18.43 |   1   |     1.50 |   1   |     7.80
            sts729 |    18.44 |   1   |     2.53 |   1   |     1.19
          circ10-3 |    18.44 |   1   |     1.12 |   1   |    14.33
        atlanta-ip |    18.44 |   0   |    13.35 |   1   |    17.24
             map20 |    18.45 |   0   |    19.16 |   1   |     8.95
             map14 |    18.45 |   0   |    18.57 |   1   |    12.04
             map18 |    18.45 |   0   |    19.42 |   1   |     9.71
          uc-case3 |    18.45 |   1   |     1.84 |   1   |     3.54
             map10 |    18.45 |   0   |    19.66 |   1   |    11.79
             map06 |    18.45 |   0   |    18.64 |   1   |    11.85
    satellites2-60 |    18.46 |   1   |     8.63 |   1   |     6.61
       neos-520729 |    18.47 |   1   |     1.44 |   1   |     2.21
       neos-957389 |    18.51 |   1   |     0.77 |   1   |     1.08
             nsr8k |    18.73 |   0   |     5.41 |   1   |     5.31
       neos-885524 |    18.75 |   1   |     1.31 |   1   |     1.27
 satellites3-40-fs |    18.80 |   0   |    39.26 |   1   |    11.04
        rocII-7-11 |    18.83 |   1   |     1.37 |   1   |     1.41
           vpphard |    18.85 |   1   |     2.23 |   1   |     6.22
             t1717 |    18.85 |   1   |     2.01 |   1   |     2.49
               ex9 |    18.93 |   1   |     1.84 |   1   |    19.90
               van |    18.99 |   1   |     2.01 |   1   |     5.49
              dc1l |    18.99 |   1   |     2.63 |   1   |     2.61
       opm2-z10-s2 |    19.05 |   1   |    33.31 |   1   |   390.08
          triptim1 |    19.08 |   0   |    29.71 |   1   |     7.69
       neos-506428 |    19.09 |   1   |     1.84 |   1   |     4.35
       neos-932816 |    19.09 |   1   |     2.50 |   1   |     2.75
          triptim2 |    19.09 |   1   |     4.32 |   1   |     7.49
          triptim3 |    19.10 |   1   |     4.14 |   1   |     6.70
         ns1116954 |    19.11 |   1   |    21.46 |   1   |   472.68
           rail507 |    19.18 |   1   |     2.97 |   1   |     2.49
         ns1904248 |    19.18 |   1   |     2.09 |   1   |     9.54
         ns1111636 |    19.19 |   1   |     2.10 |   1   |     1.83
        rocII-9-11 |    19.20 |   1   |     1.94 |   1   |     1.99
             n15-3 |    19.23 |   1   |     2.71 |   1   |     2.94
            rail01 |    19.26 |   1   |     8.13 |   1   |     9.57
        gmut-75-50 |    19.30 |   1   |     2.58 |   1   |     2.45
   pb-simp-nonunif |    19.41 |   1   |     1.91 |   1   |     4.06
       opm2-z11-s8 |    19.52 |   1   |    45.95 |   1   |   393.77
       neos-941313 |    19.67 |   1   |     4.94 |   1   |     2.60
    satellites3-40 |    19.73 |   0   |   221.62 |   1   |    20.07
      neos-1140050 |    19.76 |   1   |     9.59 |   1   |   150.90
       neos-859770 |    19.76 |   1   |     0.88 |   1   |     1.03
      netdiversion |    19.89 |   1   |     5.33 |   1   |     8.88
           rmine14 |    19.92 |   1   |    19.97 |   1   |   179.57
         momentum3 |    19.98 |   1   |     4.30 |   1   |    89.99
    buildingenergy |    19.99 |   1   |     9.61 |   1   |    17.10
           rvb-sub |    20.00 |   1   |     1.87 |   1   |     3.23
      opm2-z12-s14 |    20.02 |   1   |    90.17 |   0   |   904.86
       opm2-z12-s7 |    20.02 |   0   |   168.88 |   0   |   903.73
          vpphard2 |    20.03 |   1   |     5.10 |   1   |    24.21
             stp3d |    20.05 |   1   |    14.57 |   1   |    25.07
         eilA101-2 |    20.06 |   1   |     2.59 |   1   |     3.49
            npmv07 |    20.07 |   0   |    48.03 |   1   |    12.66
           sing245 |    20.10 |   1   |     8.81 |   1   |   106.63
         ns2118727 |    20.10 |   1   |    11.31 |   1   |    32.39
         ns2137859 |    20.11 |   1   |     7.75 |   1   |     5.27
              ex10 |    20.17 |   1   |     4.31 |   1   |    59.36
         ns1854840 |    20.27 |   1   |     6.67 |   1   |    13.37
            rail02 |    20.31 |   1   |    19.56 |   1   |    33.68
       neos-631710 |    20.35 |   1   |     4.26 |   1   |     7.58
           datt256 |    20.53 |   1   |     8.29 |   1   |   709.67
         ns1758913 |    20.89 |   1   |    12.05 |   1   |   872.01
      neos-1429212 |    20.93 |   1   |     4.08 |   1   |     9.32
  wnq-n100-mw99-14 |    20.94 |   1   |    25.63 |   1   |   775.53
         ns1853823 |    20.94 |   1   |    18.23 |   1   |    85.47
            co-100 |    21.00 |   1   |     3.55 |   1   |     4.23
            rail03 |    21.63 |   1   |    51.85 |   1   |   102.11
           n3seq24 |    21.73 |   1   |     7.41 |   1   |    16.68
       neos-476283 |    21.92 |   1   |    93.77 |   1   |    32.56
              bab3 |    21.97 |   1   |    23.74 |   1   |    24.50
           rmine21 |    22.33 |   1   |   265.36 |   0   |   919.34
         ivu06-big |    24.72 |   1   |   353.88 |   1   |   148.94
            mspp16 |    24.74 |   1   |    59.55 |   1   |   559.65
--------------------------------------------------------------------
\end{verbatim}

\begin{verbatim}
--------------------------------------------------------------------
              opf benchmark results (tol = 0.0001)
--------------------------------------------------------------------
           problem | log2(nnz)|      MadNLP      |      Ipopt        
                   |          | solved|     time | solved|     time     
--------------------------------------------------------------------
        case3_lmbd |     7.93 |   1   |     0.12 |   1   |     0.01
         case5_pjm |     8.91 |   1   |     0.18 |   1   |     0.01
       case14_ieee |    10.60 |   1   |     0.09 |   1   |     0.01
   case24_ieee_rts |    11.55 |   1   |     0.15 |   1   |     0.02
         case30_as |    11.63 |   1   |     0.10 |   1   |     0.02
       case30_ieee |    11.63 |   1   |     0.11 |   1   |     0.03
       case39_epri |    11.81 |   1   |     0.25 |   1   |     0.03
       case57_ieee |    12.59 |   1   |     0.14 |   1   |     0.02
          case60_c |    12.74 |   1   |     0.20 |   1   |     0.04
   case73_ieee_rts |    13.21 |   1   |     0.16 |   1   |     0.04
      case118_ieee |    13.81 |   1   |     0.16 |   1   |     0.05
     case89_pegase |    13.96 |   1   |     0.22 |   1   |     0.06
     case200_activ |    14.22 |   1   |     0.12 |   1   |     0.03
       case179_goc |    14.31 |   1   |     0.23 |   1   |     0.10
  case162_ieee_dtc |    14.41 |   1   |     0.36 |   1   |     0.08
      case197_snem |    14.43 |   1   |     0.13 |   1   |     0.04
      case300_ieee |    14.96 |   1   |     0.43 |   1   |     0.11
     case240_pserc |    15.08 |   1   |     1.52 |   1   |     0.63
      case588_sdet |    15.70 |   1   |     0.28 |   1   |     0.15
       case500_goc |    15.78 |   1   |     0.34 |   1   |     0.21
       case793_goc |    16.12 |   1   |     0.31 |   1   |     0.23
   case1354_pegase |    17.23 |   1   |     0.41 |   1   |     0.63
      case1888_rte |    17.58 |   1   |     2.92 |   1   |     1.41
      case1951_rte |    17.62 |   1   |     1.07 |   1   |     1.40
     case1803_snem |    17.71 |   1   |     0.48 |   1   |     1.21
      case2383wp_k |    17.78 |   1   |     0.57 |   1   |     1.28
      case2312_goc |    17.83 |   1   |     0.46 |   1   |     1.02
     case2737sop_k |    17.95 |   1   |     0.46 |   1   |     0.82
      case2736sp_k |    17.95 |   1   |     0.50 |   1   |     1.12
      case2746wp_k |    17.96 |   1   |     0.44 |   1   |     0.99
     case2746wop_k |    17.97 |   1   |     0.40 |   1   |     0.86
      case2000_goc |    18.08 |   1   |     0.44 |   1   |     1.24
      case3012wp_k |    18.08 |   1   |     0.65 |   1   |     1.83
      case3120sp_k |    18.13 |   1   |     0.59 |   1   |     1.67
      case2848_rte |    18.16 |   1   |     1.49 |   1   |     2.45
      case2868_rte |    18.17 |   1   |     2.09 |   1   |     2.81
     case2853_sdet |    18.21 |   1   |     0.60 |   1   |     2.11
      case3022_goc |    18.28 |   1   |     0.53 |   1   |     1.67
      case3375wp_k |    18.30 |   1   |     0.63 |   1   |     2.03
   case2869_pegase |    18.43 |   1   |     0.70 |   1   |     2.04
      case2742_goc |    18.45 |   1   |     2.01 |   1   |     5.59
     case4661_sdet |    18.83 |   1   |     1.01 |   1   |     3.76
      case3970_goc |    18.96 |   1   |     0.70 |   1   |     4.70
      case4917_goc |    18.99 |   1   |     0.69 |   1   |     3.42
      case4020_goc |    19.03 |   1   |     1.08 |   1   |     6.51
      case4601_goc |    19.08 |   1   |     1.01 |   1   |     5.92
      case4837_goc |    19.18 |   1   |     0.90 |   1   |     5.01
      case4619_goc |    19.25 |   1   |     0.91 |   1   |     5.89
      case6468_rte |    19.40 |   1   |    13.65 |   1   |    15.39
      case6495_rte |    19.41 |   1   |    32.78 |   1   |    15.96
      case6470_rte |    19.41 |   1   |    18.65 |   1   |     7.71
      case6515_rte |    19.41 |   1   |     5.59 |   1   |    11.69
 case5658_epigrids |    19.41 |   1   |     0.94 |   1   |     5.22
 case7336_epigrids |    19.75 |   1   |     1.05 |   1   |     6.63
     case10000_goc |    19.96 |   1   |     1.06 |   1   |    11.87
   case8387_pegase |    20.09 |   1   |     1.69 |   1   |    11.61
      case9591_goc |    20.22 |   1   |     1.82 |   1   |    20.10
   case9241_pegase |    20.23 |   1   |     1.64 |   1   |    12.14
case10192_epigrids |    20.31 |   1   |     1.60 |   1   |    15.12
     case10480_goc |    20.44 |   1   |     2.21 |   1   |    21.72
  case13659_pegase |    20.59 |   1   |     2.37 |   1   |    20.39
case20758_epigrids |    21.29 |   1   |     3.18 |   1   |    27.56
     case19402_goc |    21.34 |   1   |     4.02 |   1   |    60.80
     case30000_goc |    21.39 |   1   |     5.54 |   1   |   203.64
     case24464_goc |    21.47 |   1   |     5.09 |   1   |    40.54
case78484_epigrids |    23.20 |   1   |    16.23 |   1   |   339.30
--------------------------------------------------------------------
\end{verbatim}

\begin{verbatim}
--------------------------------------------------------------------
              opf benchmark results (tol = 1.0e-8)
--------------------------------------------------------------------
           problem | log2(nnz)|      MadNLP      |      Ipopt        
                   |          | solved|     time | solved|     time     
--------------------------------------------------------------------
        case3_lmbd |     7.93 |   1   |     0.25 |   1   |     0.01
         case5_pjm |     8.91 |   1   |     0.54 |   1   |     0.02
       case14_ieee |    10.60 |   1   |     0.13 |   1   |     0.02
   case24_ieee_rts |    11.55 |   1   |     0.38 |   1   |     0.02
       case30_ieee |    11.63 |   1   |     0.17 |   1   |     0.02
         case30_as |    11.63 |   1   |     0.12 |   1   |     0.01
       case39_epri |    11.81 |   1   |     0.43 |   1   |     0.02
       case57_ieee |    12.59 |   1   |     0.22 |   1   |     0.02
          case60_c |    12.74 |   1   |     0.30 |   1   |     0.03
   case73_ieee_rts |    13.21 |   1   |     0.39 |   1   |     0.04
      case118_ieee |    13.81 |   1   |     0.30 |   1   |     0.05
     case89_pegase |    13.96 |   1   |     0.24 |   1   |     0.07
     case200_activ |    14.22 |   1   |     0.30 |   1   |     0.05
       case179_goc |    14.31 |   1   |     0.33 |   1   |     0.13
  case162_ieee_dtc |    14.41 |   1   |     0.36 |   1   |     0.08
      case197_snem |    14.43 |   1   |     0.48 |   1   |     0.07
      case300_ieee |    14.96 |   1   |     0.73 |   1   |     0.12
     case240_pserc |    15.08 |   1   |     2.00 |   1   |     0.65
      case588_sdet |    15.70 |   1   |     0.39 |   1   |     0.19
       case500_goc |    15.78 |   1   |     0.60 |   1   |     0.23
       case793_goc |    16.12 |   1   |     0.81 |   1   |     0.28
   case1354_pegase |    17.23 |   1   |     0.83 |   1   |     0.76
      case1888_rte |    17.58 |   0   |    13.39 |   1   |     3.43
      case1951_rte |    17.62 |   1   |    15.49 |   1   |     1.71
     case1803_snem |    17.71 |   1   |     1.14 |   1   |     1.40
      case2383wp_k |    17.78 |   1   |     0.73 |   1   |     1.46
      case2312_goc |    17.83 |   1   |     0.83 |   1   |     1.22
      case2736sp_k |    17.95 |   1   |     0.57 |   1   |     1.17
     case2737sop_k |    17.95 |   1   |     0.52 |   1   |     1.04
      case2746wp_k |    17.96 |   1   |     0.60 |   1   |     1.19
     case2746wop_k |    17.97 |   1   |     0.67 |   1   |     1.05
      case3012wp_k |    18.08 |   1   |     0.91 |   1   |     1.97
      case2000_goc |    18.08 |   1   |     0.81 |   1   |     1.35
      case3120sp_k |    18.13 |   1   |     0.85 |   1   |     1.92
      case2848_rte |    18.16 |   1   |     7.06 |   1   |     2.76
      case2868_rte |    18.17 |   1   |    29.18 |   1   |     2.99
     case2853_sdet |    18.21 |   1   |     0.91 |   1   |     1.82
      case3022_goc |    18.28 |   1   |     1.15 |   1   |     2.03
      case3375wp_k |    18.30 |   1   |     1.41 |   1   |     2.37
   case2869_pegase |    18.43 |   1   |     0.94 |   1   |     2.50
      case2742_goc |    18.45 |   1   |     4.43 |   1   |     6.01
     case4661_sdet |    18.83 |   1   |     2.13 |   1   |     4.16
      case3970_goc |    18.96 |   1   |     1.47 |   1   |     5.18
      case4917_goc |    18.99 |   1   |     1.48 |   1   |     4.24
      case4020_goc |    19.03 |   1   |     1.88 |   1   |     7.04
      case4601_goc |    19.08 |   1   |     2.13 |   1   |     6.23
      case4837_goc |    19.18 |   1   |     1.46 |   1   |     5.23
      case4619_goc |    19.25 |   1   |     1.50 |   1   |     6.15
      case6468_rte |    19.40 |   1   |    11.27 |   1   |    13.55
      case6470_rte |    19.41 |   1   |    20.15 |   1   |     8.39
 case5658_epigrids |    19.41 |   1   |     1.65 |   1   |     5.90
      case6495_rte |    19.41 |   1   |    50.00 |   1   |    16.25
      case6515_rte |    19.41 |   1   |    13.51 |   1   |    12.47
 case7336_epigrids |    19.75 |   1   |     1.89 |   1   |     7.64
     case10000_goc |    19.96 |   1   |     2.39 |   1   |    13.54
   case8387_pegase |    20.09 |   1   |     3.28 |   1   |    13.16
      case9591_goc |    20.22 |   1   |     3.21 |   1   |    22.01
   case9241_pegase |    20.23 |   0   |   114.25 |   1   |    13.93
case10192_epigrids |    20.31 |   1   |     2.84 |   1   |    17.14
     case10480_goc |    20.44 |   1   |     3.24 |   1   |    23.52
  case13659_pegase |    20.59 |   1   |     3.41 |   1   |    17.56
case20758_epigrids |    21.29 |   1   |     7.98 |   1   |    31.55
     case19402_goc |    21.34 |   1   |     5.27 |   1   |    62.56
     case30000_goc |    21.39 |   1   |     9.59 |   1   |    98.95
     case24464_goc |    21.47 |   1   |     4.96 |   1   |    43.97
case78484_epigrids |    23.20 |   1   |    19.62 |   1   |   365.90
--------------------------------------------------------------------
\end{verbatim}

\begin{verbatim}
--------------------------------------------------------------------
              cops benchmark results (tol = 0.0001)
--------------------------------------------------------------------
           problem | log2(nnz)|      MadNLP      |      Ipopt        
                   |          | solved|     time | solved|     time     
--------------------------------------------------------------------
     camshape-1600 |    14.10 |   1   |     0.19 |   1   |     0.66
     camshape-3200 |    15.10 |   1   |     0.23 |   1   |     3.56
         robot-400 |    15.24 |   1   |     0.37 |   1   |     8.34
     camshape-6400 |    16.10 |   1   |     0.37 |   1   |    16.63
        marine-400 |    16.10 |   1   |     0.26 |   1   |     0.26
         robot-800 |    16.24 |   1   |     3.91 |   1   |     9.54
          elec-100 |    16.67 |   1   |     0.90 |   1   |     1.49
     steering-3200 |    17.10 |   1   |     0.29 |   1   |     0.45
    camshape-12800 |    17.10 |   1   |     0.27 |   1   |    66.00
        marine-800 |    17.10 |   1   |     0.29 |   1   |     0.61
        gasoil-800 |    17.22 |   1   |     0.23 |   1   |     0.46
        robot-1600 |    17.24 |   1   |     4.61 |   1   |     3.59
       rocket-3200 |    17.93 |   1   |     0.38 |   1   |     2.19
        pinene-800 |    18.01 |   1   |     0.48 |   1   |     0.97
       marine-1600 |    18.10 |   1   |     0.46 |   1   |     2.88
     steering-6400 |    18.10 |   1   |     0.62 |   1   |     1.06
    camshape-25600 |    18.10 |   1   |     0.30 |   1   |   301.15
       gasoil-1600 |    18.22 |   1   |     0.47 |   1   |     1.22
        robot-3200 |    18.24 |   0   |    19.21 |   1   |   103.29
   bearing-200,200 |    18.63 |   1   |     0.34 |   1   |     0.78
          elec-200 |    18.68 |   1   |     0.74 |   1   |     4.52
       rocket-6400 |    18.93 |   1   |     0.91 |   1   |    11.88
       pinene-1600 |    19.01 |   1   |     0.80 |   1   |     2.52
       marine-3200 |    19.10 |   1   |     1.05 |   1   |    12.61
    steering-12800 |    19.10 |   1   |     0.97 |   1   |     3.66
       gasoil-3200 |    19.22 |   1   |     0.99 |   1   |     6.14
        robot-6400 |    19.24 |   1   |    11.48 |   0   |   177.21
   bearing-300,300 |    19.79 |   1   |     1.31 |   1   |     1.95
      rocket-12800 |    19.93 |   1   |     1.76 |   1   |    19.82
       pinene-3200 |    20.01 |   1   |     2.17 |   1   |     4.87
       marine-6400 |    20.10 |   1   |     2.07 |   1   |    59.12
    steering-25600 |    20.10 |   1   |     1.78 |   1   |     7.89
       gasoil-6400 |    20.21 |   1   |     1.63 |   1   |    18.30
   bearing-400,400 |    20.62 |   1   |     0.99 |   1   |     3.73
          elec-400 |    20.68 |   1   |     1.90 |   1   |    28.20
      rocket-25600 |    20.93 |   1   |     3.64 |   1   |   113.01
       pinene-6400 |    21.01 |   1   |     2.95 |   1   |    15.86
    steering-51200 |    21.10 |   1   |     3.52 |   1   |    20.09
      gasoil-12800 |    21.21 |   1   |     3.33 |   1   |    24.19
   bearing-600,600 |    21.79 |   1   |     2.30 |   1   |     9.47
      rocket-51200 |    21.93 |   1   |     6.19 |   1   |   856.66
      pinene-12800 |    22.01 |   1   |     5.79 |   1   |    71.92
   bearing-800,800 |    22.61 |   1   |     4.17 |   1   |    19.09
          elec-800 |    22.68 |   1   |     9.27 |   1   |   263.42
         elec-1600 |    24.68 |   1   |    14.56 |   0   |   909.28
--------------------------------------------------------------------
\end{verbatim}

\begin{verbatim}
--------------------------------------------------------------------
              cops benchmark results (tol = 1.0e-8)
--------------------------------------------------------------------
           problem | log2(nnz)|      MadNLP      |      Ipopt        
                   |          | solved|     time | solved|     time     
--------------------------------------------------------------------
     camshape-1600 |    14.10 |   1   |     0.36 |   1   |     1.13
     camshape-3200 |    15.10 |   1   |     0.45 |   1   |     3.89
         robot-400 |    15.24 |   1   |     1.05 |   1   |     1.00
        marine-400 |    16.10 |   1   |     0.30 |   1   |     0.28
     camshape-6400 |    16.10 |   1   |     1.07 |   1   |    17.16
         robot-800 |    16.24 |   1   |     0.74 |   1   |     2.41
          elec-100 |    16.67 |   1   |     0.53 |   1   |     1.50
     steering-3200 |    17.10 |   1   |     0.53 |   1   |     0.52
    camshape-12800 |    17.10 |   1   |     1.20 |   1   |    95.98
        marine-800 |    17.10 |   1   |     0.36 |   1   |     0.81
        gasoil-800 |    17.22 |   1   |     0.55 |   1   |     0.53
        robot-1600 |    17.24 |   1   |     1.73 |   1   |    33.31
       rocket-3200 |    17.93 |   1   |     2.48 |   1   |     1.53
        pinene-800 |    18.01 |   1   |     0.50 |   1   |     1.09
       marine-1600 |    18.10 |   1   |     0.58 |   1   |     3.37
     steering-6400 |    18.10 |   1   |     0.99 |   1   |     1.32
    camshape-25600 |    18.10 |   1   |     1.02 |   1   |   451.05
       gasoil-1600 |    18.22 |   1   |     0.83 |   1   |     1.26
        robot-3200 |    18.24 |   1   |     2.22 |   1   |   147.12
   bearing-200,200 |    18.63 |   1   |     0.39 |   1   |     1.23
          elec-200 |    18.68 |   1   |     2.14 |   1   |     4.58
       rocket-6400 |    18.93 |   1   |     1.50 |   1   |     3.63
       pinene-1600 |    19.01 |   1   |     0.67 |   1   |     2.66
       marine-3200 |    19.10 |   1   |     1.26 |   1   |    23.98
    steering-12800 |    19.10 |   1   |     3.06 |   1   |     4.04
       gasoil-3200 |    19.22 |   1   |     1.66 |   1   |     5.64
        robot-6400 |    19.24 |   1   |     5.72 |   0   |   900.32
   bearing-300,300 |    19.79 |   1   |     0.63 |   1   |     3.17
      rocket-12800 |    19.93 |   1   |     2.81 |   1   |    24.24
       pinene-3200 |    20.01 |   1   |     1.47 |   1   |     5.18
       marine-6400 |    20.10 |   1   |     3.39 |   1   |   123.19
    steering-25600 |    20.10 |   1   |     2.99 |   1   |     8.88
       gasoil-6400 |    20.21 |   1   |     7.75 |   1   |    11.23
   bearing-400,400 |    20.62 |   1   |     1.12 |   1   |     5.82
          elec-400 |    20.68 |   1   |     2.99 |   1   |    28.05
      rocket-25600 |    20.93 |   1   |    19.02 |   1   |    33.88
       pinene-6400 |    21.01 |   1   |     3.64 |   1   |    19.28
    steering-51200 |    21.10 |   1   |     6.69 |   1   |    23.00
      gasoil-12800 |    21.21 |   1   |    24.91 |   1   |    28.71
   bearing-600,600 |    21.79 |   1   |     2.40 |   1   |    14.62
      rocket-51200 |    21.93 |   1   |    27.43 |   1   |    81.64
      pinene-12800 |    22.01 |   1   |     6.35 |   1   |    81.15
   bearing-800,800 |    22.61 |   1   |     4.30 |   1   |    29.51
          elec-800 |    22.68 |   1   |    10.56 |   1   |   329.08
         elec-1600 |    24.68 |   1   |    53.04 |   0   |   909.32
--------------------------------------------------------------------
\end{verbatim}

\end{document}